\documentclass[reqno,11pt]{amsart}
\usepackage{latexsym}
\usepackage{amsfonts}
\usepackage{mathrsfs}
\usepackage{graphicx}

\newtheorem{lem}{Lemma}[section]
\newtheorem{thm}{Theorem}[section]

\newtheorem{rem}{Remark}[section]
\newtheorem{exa}{Example}[section] 

\newcommand{\bz}{{ \bf 0 }}

\newcommand{\Rs}{\mathbb{R}}

\newcommand{\PP}{{\mathcal{P}} }

\newcommand{\rr}{\bar{r} }

\newcommand{\bpr}{{\bf Proof.} \hspace{1 em}}
\newcommand{\epr}{ \\ \hspace*{4.5in} $\Box$ }
\newcommand{\beq}{ \begin{equation} }
\newcommand{\eeq}{ \end{equation} }

\numberwithin{equation}{section}

\begin{document}

\bibliographystyle{plain}
\title[On Stress matrices of Chordal Bar Frameworks]
{On Stress matrices of Chordal Bar Frameworks in General Position}

\author{ A. Y. Alfakih} 

\address{Department of Mathematics and Statistics \\
          University of Windsor \\
          Windsor, Ontario N9B 3P4 \\
          Canada. } 
\email{alfakih@uwindsor.ca} 
  
\thanks{Research supported by the Natural Sciences and Engineering
         Research Council of Canada.} 
          
\date{\today}  

\subjclass[2000]{05C50, 52C25, 15B99}  

\keywords {Bar frameworks, chordal graphs, stress matrices, 
Gale matrices, universal rigidity, global rigidity and points in
general position.} 

\date{\today}  

\begin{abstract}
A bar framework in $\Rs^r$, denoted by $G(p)$, 
is a simple connected graph $G$ whose vertices 
are points  $p^1, \ldots, p^n$ in $\Rs^r$ that affinely span $\Rs^r$,
and whose edges are line segments between pairs of these points.  
In this paper, we use stress matrices to characterize the universal and global
rigidities of chordal bar frameworks in general position in $\Rs^r$, i.e., 
bar frameworks where graph $G$ is chordal and the  
points $p^1, \ldots, p^n$ are in general position in $\Rs^r$.  
We also prove that if a chordal bar framework in $\Rs^r$ admits a stress
matrix of rank $n-r-1$ with generic rank profile, then it admits 
a positive semidefinite stress matrix of rank $n-r-1$.  
\end{abstract}

\maketitle

\section{Introduction}
\label{int}

A {\em bar framework} (a framework for short) in $r$-dimensional
Euclidean space, 
denoted by $G(p)$, is a connected simple graph $G=(V,E)$
(a graph with no loops or multiple edges)
whose vertices are points $p^1, \ldots, p^n$ in $\Rs^r$ that affinely span $\Rs^r$,
and whose edges are line segments between pairs of these points.    
The points $p^1, \ldots, p^n$, which collectively are denoted by $p$, 
are referred to as  the {\em configuration} of the framework.
Throughout this paper, $V(G)$ and $E(G)$ denote the vertex set and the edge set of
graph $G$ respectively; and  
\bz  denotes the zero matrix or the zero vector of the
appropriate dimension. 

We say that frameworks $G(p)$ and $G(q)$ in $\Rs^r$ are 
{\em congruent} if $||q^i-q^j||$= $||p^i-p^j||$ for all
$i,j=1,\ldots,n$, where $||.||$ denotes the Euclidean norm. 
For example, if configuration $q$ is obtained from configuration $p$
by a rigid motion, then $G(p)$ and $G(q)$ are congruent. On the other hand, 
frameworks $G(p)$ in $\Rs^r$ and $G(q)$ in $\Rs^s$ are said to be  
{\em equivalent} if the corresponding edges in the two frameworks have 
the same length, i.e., if $||q^i-q^j||$= $||p^i-p^j||$ for all
$(i,j) \in E(G)$. 

A framework $G(p)$ in $\Rs^r$ is said to be {\em globally rigid}
if every framework $G(q)$ in the same dimension which  
is equivalent to $G(p)$, is in fact congruent to $G(p)$. Furthermore,  
if every framework $G(q)$ in any dimension which   
is equivalent to $G(p)$ is congruent to $G(p)$, then   
$G(p)$ is said to be {\em universally rigid}. 
Obviously, universal rigidity implies global rigidity. 

A real valued function $\omega$ on the edge set $E(G)$ such that 
\beq \label{defw} 
\sum_{j : (i,j) \in E(G)} \omega_{ij}(p^i-p^j) = \bz \mbox{ for all } i=1,\ldots,n
\eeq
is called an {\em equilibrium stress} of framework $G(p)$. Furthermore,
given an equilibrium stress $\omega$, 
the $n \times n$ symmetric matrix $S=(s_{ij})$ where
\beq \label{defS} 
 s_{ij} = \left\{ \begin{array}{ll} -\omega_{ij} & \mbox{if } (i,j) \in E(G),\\
                       0   & \mbox{if $i \neq j$ and } (i,j) \not \in E(G), \\
                \sum_{k: (i,k) \in E(G)} \omega_{ik} & \mbox{if } i=j,
                     \end{array} \right. 
\eeq
is called a {\em stress matrix} of $G(p)$.

The problems of universal and global rigidities of bar frameworks
have many important applications in molecular conformations and sensor networks. 
As a result, they have received a great deal of attention. Furthermore,
under the assumption of generic configurations, 
global and universal rigidities have nice characterizations in terms of
stress matrices \cite{con99, hen92,alf10,gt09, con05,ght10}.
A configuration $p$ is said to be generic if the coordinates of
$p^1, \ldots,p^n$ are algebraically independent over the integers.

The assumption of a generic configuration is quite strong. A study
of the problem of universal rigidity under the weaker assumption
of points in general position was initiated in \cite{ay10, aty10c}.   
We say that a configuration $p$ in $\Rs^r$ is in {\em general position} if
every ($r+1$) of the points $p^1, \ldots, p^n$ are affinely independent. 
Recall that points $p^1, \ldots, p^{r+1}$ are {\em affinely independent} if the 
trivial solution $\lambda_1=\cdots=\lambda_{r+1}=0$ is the only   
solution for the system of equations 
\beq \label{defai} 
\sum_{i=1}^{r+1} \lambda_i \left[ \begin{array}{c} p^i \\ 1  \end{array} \right] = \bz.  
\eeq 
For example, points in $\Rs^2$ are in general position if no three of them are
collinear. An $n \times n$ real symmetric matrix $A$ is {\em positive semidefinite} if  
$x^T A x \geq 0$ for all $x \in \Rs^n$. 

The following two theorems are the main results in \cite{ay10} and \cite{aty10c}.   

\begin{thm}[Alfakih and Ye \cite{ay10}] \label{suffgp} 
Let $G(p)$ be a bar framework on $n$ vertices 
in general position in $\Rs^r$. Then $G(p)$ is universally rigid if
$G(p)$ admits a positive semidefinite stress matrix $S$ of rank $n-r-1$. 
\end{thm}

\begin{thm}[Alfakih et al \cite{aty10c}] \label{nectri} 
Let $G(p)$ be a bar framework on $n$ vertices 
in general position in $\Rs^r$, where $G$ is  
an ($r+1$)-lateration graph. Then 
$G(p)$ admits a positive semidefinite stress matrix $S$ of rank $n-r-1$. 
\end{thm}

This paper builds on some of the ideas in \cite{aty10c} to characterize,
the universal and global rigidities of chordal frameworks in general position,   
i.e., frameworks whose graphs are chordal, and whose configurations 
are in general position. 
A graph $G$ is said to be {\em chordal} (also called {\em triangulated, rigid circuit}) 
if every cycle of $G$ of length $\geq 4$ has a chord, i.e.,   
an edge that connects two non-consecutive vertices on the cycle.

A graph $G$ is said to be {\em $k$-vertex connected} 
iff $G$ is the complete graph on $k+1$ vertices, 
or $|V(G)| \geq k+2$ and the deletion of any $k-1$ vertices leaves
$G$ connected.  
The next two theorems are the main results of this paper. 

\begin{thm} \label{main1} 
Let  $G(p)$ be a chordal bar framework on $n$ vertices 
in general position in $\Rs^r$. Then the following statements are
equivalent:
\begin{enumerate}
\item $G$ is ($r+1$)-vertex connected.
\item $G(p)$ admits a positive semidefinite stress matrix $S$ of  
rank $n-r-1$. 
\item $G(p)$ is universally rigid.
\item $G(p)$ is globally rigid. 
\end{enumerate}
\end{thm}

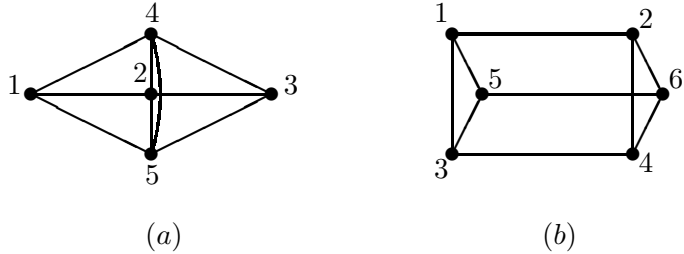
\begin{figure}[t] 
\thicklines 
\setlength{\unitlength}{0.8mm}
\begin{picture}(40,75)(-20,-40) 
\put(-60,0){\circle*{2}} 
\put(-40,0){\circle*{2}} 
\put(-20,0){\circle*{2}} 
\put(-40,10){\circle*{2}} 
\put(-40,-10){\circle*{2}} 

\put(-60,0){\line(1,0){20}}
\put(-60,0){\line(2,1){20}}
\put(-60,0){\line(2,-1){20}}
\put(-40,0){\line(1,0){20}}
\put(-40,0){\line(0,1){10}}
\put(-40,0){\line(0,-1){10}}
\put(-40,10){\line(2,-1){20}}
\put(-40,-10){\line(2,1){20}}
\qbezier(-40,-10)(-37,0) (-40,10)

\put(-64,0){$1$}
\put(-43,2){$2$}
\put(-18,0){$3$}
\put(-41,12){$4$}
\put(-41,-15){$5$}

\put(15,0){\circle*{2}} 
\put(45,0){\circle*{2}} 
\put(10,10){\circle*{2}} 
\put(10,-10){\circle*{2}} 
\put(40,10){\circle*{2}} 
\put(40,-10){\circle*{2}} 

\put(15,0){\line(1,0){30}}
\put(10,10){\line(0,-1){20}}
\put(10,10){\line(1,0){30}}
\put(10,10){\line(1,-2){5}}
\put(40,10){\line(0,-1){20}}
\put(40,10){\line(1,-2){5}}
\put(40,-10){\line(1,2){5}}
\put(10,-10){\line(1,2){5}}
\put(10,-10){\line(1,0){30}}

\put(7,12){$1$}
\put(41,11){$2$}
\put(7,-14){$3$}
\put(41,-13){$4$}
\put(16,1){$5$}
\put(46,1){$6$}

\put(-41,-25){$(a)$}
\put(25,-25){$(b)$}
\end{picture}
\caption{Two non-globally rigid frameworks in $\Rs^2$. 
     The edge $(4,5)$ in framework $(a)$ is shown as an arc to make
edges $(2,4)$ and $(2,5)$ visible. } 
\label{f1}
\end{figure}

The proof of Theorem \ref{main1} is given in Section \ref{sec3}. 
We remark here that the equivalence of statements $(1)$ and $(3)$ in Theorem \ref{main1}
also follows from a result by Bakonyi and Johnson concerning Euclidean distance
matrix completions \cite{bj95}. 
The assumptions of chordality and general position cannot
be dropped in Theorem \ref{main1}. 
Figure \ref{f1} depicts two frameworks in $\Rs^2$ that are not globally rigid. 
Framework $(a)$ is 3-vertex connected and chordal, but it is not in general position;
while framework $(b)$ is in general position and 3-vertex connected, but it is  
not chordal. 

The {\em $k$th leading principal minor} of an $n \times n$ matrix $A$ is
the determinant of the square submatrix obtained by deleting the last $(n-k)$ rows
and columns of $A$. Let $A$ be a given symmetric matrix of rank $k$, 
then $A$ is said to have
{\em generic rank profile} if its first $k$ leading principal minors
are nonzero. 

\begin{thm} \label{main2} 
Let  $G(p)$ be a chordal bar framework on $n$ vertices 
in general position in $\Rs^r$. If 
$G(p)$ admits a stress matrix $S$ of rank $n-r-1$ with generic rank profile, 
then it admits a positive semidefinite stress matrix $S'$ of rank $n-r-1$. 
\end{thm}

The proof of Theorem \ref{main2} is given in Section \ref{sec4}. 

\section{Preliminaries}
\label{pre}
In this section we present some necessary mathematical preliminaries.   
In particular we review some facts concerning chordal graphs and Gale matrices, 
which will be used in our proofs. 

\subsection{Chordal Graphs}

A simple graph $G=(V,E)$ is {\em complete} if the nodes of $G$
are pairwise adjacent. A  complete subgraph is called  a {\em clique}. 
For a vertex $v$ in $G$, $N(v)$ denotes the set of nodes of $G$ that
are adjacent to $v$. A vertex $v$ of $G$ is said to be {\em simplicial} if 
$N(v)$ is a clique in $G$. A permutation of the vertices of $G$, 
$\pi(1), \pi(2), \ldots, \pi(n)$, is called a 
{\em perfect elimination ordering (PEO)} of $G$ if for each $i=1,\ldots, n-1$, 
$N(\pi(i))$ is simplicial in the graph induced by 
$\{\pi(i),\pi({i+1}), \ldots, \pi(n) \}$. 

Among the many well-known characterizations of chordal graphs 
\cite{bp93,gol04}, the one presented in the following lemma 
is the most useful for the purposes of this paper. 

\begin{thm}[Fulkerson and Gross \cite{fg65}]
A graph $G$ is chordal if and only if $G$ has a perfect elimination ordering. 
\end{thm}

A PEO of a chordal graph $G$ on $n$ vertices can be found by using, for example, the 
{\em Maximum Cardinality Search (MCS)} algorithm \cite{tar76,ty84}. 
This algorithm runs in 
$O(n+|E(G)|)$ time, and obtains a PEO of $G$  as follows. 
Initially all vertices of $G$ are unlabeled.
Arbitrarily select a vertex and label it $\pi(n)$. Then for $i=n-1, n-2, \ldots,1$,
select a vertex which is adjacent to the largest number of already 
labeled vertices, breaking ties arbitrarily, and label it $\pi(i)$. 
One can then prove that $\pi(1), \ldots, \pi(n)$ is a PEO of $G$.
For ease of notation and without loss of generality, we assume 
throughout the paper that $1,2,\ldots,n$ is a PEO of $G$. 

The following lemma will be needed in the sequel.  

\begin{lem}[Lin et al \cite{lcc07}] \label{lin} 
Let $G$ be a chordal graph on $n$ vertices and let
$1,2,\ldots,n$ be a PEO of $G$. Furthermore, let
$\hat{N}(j)=\{ i : i > j \mbox{ and } (i,j) \in E(G) \}$.  
Then $G$ is $k$-vertex connected if and only if  
$|\hat{N}(j)| \geq k$ for all $j=1,\ldots, n-k$, 
where $|\hat{N}(j)|$ denotes the cardinality of  $\hat{N}(j)$ 
\end{lem}

\subsection{Gale Matrices and Stress Matrices}

Gale matrices are intimately related to stress matrices. 
Given a framework $G(p)$ in $\Rs^r$,  
the $(r+1) \times n$ matrix 
\beq \label{defPe} 
{\PP} := 
 \left[ \begin{array}{ccc} p^1 &  \ldots & p^n \\ 
                           1   &  \ldots & 1 
         \end{array}  \right], 
\eeq 
is called the {\em extended configuration matrix} of $G(p)$. 
$\PP$ has full row rank since $p^1,\ldots,p^n$ affinely span $\Rs^r$. 
Furthermore, under the general position assumption, every 
$(r+1) \times (r+1)$ submatrix of $\PP$ is nonsingular. 
Let $\rr$ be the nullity of $\PP$, i.e., the dimension of its null space.
Then 
\beq
\rr = n- r -1. 
\eeq
For $\rr \geq 1$, any $n \times \rr$
matrix $Z$ whose columns form a basis of the null space of $\PP$ 
is called a {\em Gale matrix} of $G(p)$ \cite{gal56}. 
Gale matrices (or Gale transform) are widely used in polytope theory \cite{gru67}.
Two remarks are in order here. First, 
the columns of $Z$ express the affine dependencies among 
the points $p^1,\ldots,p^n$. Second, if $Z$ is a Gale matrix
of $G(p)$ and $Q$ is any nonsingular $\rr \times \rr$ matrix.
Then $Z'=ZQ$ is also a Gale matrix of $G(p)$. This fact will be used 
very often in the sequel. 

It is not hard to see from (\ref{defS}) that an $n \times n $ symmetric
matrix $S=(s_{ij})$ is a stress matrix of $G(p)$ if and only if   
\begin{enumerate}
\item $  \PP S = \bz $, and 
\item $  s_{ij}= 0$  for all $ij: i \neq j, (i,j) \not \in E(G)$, 
\end{enumerate}
where $\PP$ is the  extended configuration matrix of $G(p)$. 
The following lemma establishes the relationship between
Gale matrices and stress matrices. 

\begin{lem}[Alfakih \cite{alf10}] \label{SZ}
Let $G(p)$ be a bar framework on $n$ vertices in $\Rs^r$, and   
let $S$ and $Z$  be, respectively, a stress matrix  and a Gale matrix 
of $G(p)$. Then $S= Z \Psi Z^T $ for some
$\rr \times \rr $ symmetric  matrix $\Psi$. Conversely, let 
$z^i$ denote the $i$th row of $Z$ and let $\Psi'$ be a symmetric
matrix such that  $(z^i)^T \Psi' z^j =0$ for all 
$ij, i \neq j, (i,j) \not \in E(G)$. Then $S'= Z \Psi' Z^T$ is a stress matrix
of $G(p)$. 
\end{lem}

It immediately follows from Lemma \ref{SZ} that rank 
$S \leq \rr$, and rank $S = \rr$ if and only if
$\Psi$ is nonsingular.  
 
The next lemma shows that Gale matrices have a useful property under
the general position assumption.

\begin{lem} \label{Zgp} 
Let $G(p)$ be a framework on $n$ vertices in general position in $\Rs^r$, 
and let $Z$ be a Gale matrix of $G(p)$. Then every $\rr \times \rr$
submatrix of $Z$ is nonsingular. 
\end{lem}

For a proof of Lemma \ref{Zgp} see, for example, \cite{alf07a}.
The next theorem is an immediate corollary of Lemma \ref{Zgp}.

\begin{thm}
Let $G(p)$ be a framework on $n$ vertices in general position in $\Rs^r$ 
and let $S$ be a stress matrix of $G(p)$ of rank $\rr$. Then 
every $\rr \times \rr$ submatrix of $S$ is nonsingular.
\end{thm}

\bpr
Let $Z$ be a Gale matrix of $G(p)$. Then $S=Z \Psi Z^T$ for some $\rr \times \rr$
symmetric matrix $\Psi$. Moreover, $\Psi$ is nonsingular since rank $S = \rr$.
Let $\alpha$ and $\beta$ be two subsets of $ \{1,2,\ldots,n\}$ of cardinality
$\rr$, and let $S_{\alpha,\beta}$ be the submatrix of $S$ whose rows and columns
are indexed by $\alpha$ and $\beta$ respectively. Further, let $Z_{\alpha}$ be the
submatrix of $Z$ of order $\rr$ whose rows are indexed by $\alpha$.  
Then
\[
S_{\alpha,\beta} = Z_{\alpha} \Psi Z^T_{\beta}.
\]

Thus it follows from Lemma \ref{Zgp} that $Z_{\alpha}$ and $Z^T_{\beta}$ 
are nonsingular and hence $S_{\alpha,\beta}$ is nonsingular.
\epr

\subsection{Gale Matrices and Chordal Frameworks}

Let $G(p)$ be a chordal framework and let
$1,2,\ldots,n$ be a PEO of $G$. 
As we saw earlier, $G(p)$ does not admit a unique Gale matrix.
Of particular interest to our purposes 
is a Gale matrix $Z=(z_{ij})$ that satisfies the following
property:

\begin{eqnarray*} 
(A) && z_{ij} =  \left\{ \begin{array}{ll} 1 & \mbox{if } i=j, \\  
                                    0 & \mbox{if } i < j , \\
                      0 & \mbox{if } i> j \mbox{ and } (i,j) \not \in E(G). 
  \end{array} \right.
\end{eqnarray*} 
That is, $Z$ is of the form
\[
\left[ \begin{array}{ccccc} 1 & & & & \\ * & 1 & & & \\ 
                            & * & 1 &  & \\
                            * &* &* &\ddots & \\
                            * & & * & & 1\\
                             & & * &* & *\\
                            * &* &  &* & \\
                            \vdots & \vdots & \vdots &\vdots & \vdots \\
                             &* & *  &* &* \\
       \end{array} \right],
\]
where zero entries are left blank and possibly nonzero entries, i.e., 
entries $ij$ where $i>j$ and $(i,j) \in E(G)$, are 
indicated with (*). 

The importance of Property (A) is illustrated in the next Lemma.  

\begin{lem} \label{Astress} 
Let $G(p)$ be a chordal framework on $n$ vertices in $\Rs^r$,
and assume that $1,2,\ldots,n$ is a PEO of $G$. If 
$G(p)$ admits a Gale matrix $Z$ that satisfies Property (A), then  
$G(p)$ admits  a positive
semidefinite stress matrix $S$ with rank $\rr$.  
\end{lem}

\bpr
Let $Z$ be a Gale matrix of $G(p)$ that satisfies Property (A) and
let $S=ZZ^T$. Then obviously $S$ is symmetric, positive semidefinite, and
of rank $\rr$. Thus, it suffices to show that 
$s_{ij} =0$ for all $ij: i \neq j, (i,j) \not \in E(G)$. 
To this end, let $(i,j) \not \in E(G)$ and $i < j$. Then    
\[ s_{ij} = \sum_{k=1}^{\rr} z_{ik} z_{jk}  = 
\sum_{k= 1}^{ \min(i,\rr)} z_{ik} z_{jk} 
\]
since $z_{ik}=0$ for $k > i$.  
But for any $k : 1 \leq k < i < j$ and $(i,j) \not \in E(G)$, either   
$i \not \in \hat{N}(k)$ (which implies that $z_{ik}=0$ ) or $j \not \in \hat{N}(k)$ 
(which implies that $z_{jk}=0$) since $\hat{N}(k)$ is
a clique. Therefore, 
\[ 
\sum_{k= 1}^{ \min(i,\rr)} z_{ik} z_{jk} = 
              \left\{ \begin{array}{ll} 0 & \mbox{if } \rr <  i , \\
                           z_{ii} z_{ji} & \mbox{if } \rr \geq i. 
                      \end{array} \right.
\]
But $z_{ji}=0$ since $(j,i) \not \in E(G)$ and $i < j$. Therefore, 
$s_{ij} =0$ for all $ij: i < j, (i,j) \not \in E(G)$ and 
the result follows since $S$ is symmetric.   
\epr 

\section{Proof of Theorem \ref{main1} }
\label{sec3}

The following lemma is needed for our proof.

\begin{lem} \label{conA} 
Let $G(p)$ be a chordal framework on $n$ vertices  in general 
position in $\Rs^r$ and assume that $1,2,\ldots,n$ is a PEO of $G$. 
If $G$ is ($r+1$)-vertex connected then $G(p)$ admits a Gale matrix $Z$ 
that satisfies Property (A).   
\end{lem}
\bpr
Recall that $\hat{N}(j)=\{ i : i > j \mbox{ and } (i,j) \in E(G) \}$.  
Since $G$ is ($r+1$)-vertex connected, it follows from Lemma \ref{lin} that 
$|\hat{N}(j)| \geq r+1$ for all $j=1, \ldots, \rr$.  
Under the general position assumption, for each $j=1, \ldots, \rr$,  
the following system of equations 
\beq \label{sysZ} 
\left[ \begin{array}{c} p^j \\ 1 \end{array} \right]+ \sum_{i \in \hat{N}(j)} x_{ij}  
 \left[ \begin{array}{c} p^i \\ 1 \end{array} \right] = \bz   
\eeq 
has a solution $\hat{x}_{ij}$. Thus, let 
$Z=(z_{ij})$ be the $n \times \rr$ matrix where 
\beq  
z_{ij} =  \left\{ \begin{array}{cl} 1 & \mbox{if } i=j, \\  
                          \hat{x}_{ij} & \mbox{if } i \in \hat{N}(j), \\  
                                    0 & \mbox{otherwise}. 
  \end{array} \right.
\eeq 
Therefore, $Z$ is a Gale matrix of $G(p)$ that satisfies Property (A). 
\epr

\noindent {\bf Proof of Theorem \ref{main1}. }

 $(1) \Rightarrow (2)$  follows from Lemmas \ref{conA} and \ref{Astress},  
 $(2) \Rightarrow (3)$  follows from Theorem \ref{suffgp},   
 $(3) \Rightarrow (4)$  is obvious,   
 $(4) \Rightarrow (1)$ follows from Hendrickson's necessary conditions
for global rigidity \cite{hen92}. If $G$ is not ($r+1$)-vertex connected,
then there exists a vertex cut $X$ of size $\leq r$ whose removal
disconnects $G(p)$ into at least two subframeworks.  
Since the vertices in $X$ are contained in an ($r-1$)-dimensional 
hyperplane, one could obtain another framework $G(q)$ that is equivalent,
but not congruent, to $G(p)$ by   
reflecting one of the subframeworks with respect to this hyperplane. 
Hence, $G(p)$ is not globally rigid.
\epr

\begin{figure}[t] 
\thicklines 
\setlength{\unitlength}{0.8mm}
\begin{picture}(40,70)(-20,-30) 
\put(-15,0){\circle*{2}} 
\put(0,0){\circle*{2}} 
\put(-15,15){\circle*{2}} 
\put(0,15){\circle*{2}} 
\put(-23,7.5){\circle*{2}} 
\put(8,7.5){\circle*{2}} 
\put(-7.5,23){\circle*{2}} 
\put(-7.5,-8){\circle*{2}} 

\put(-15,15){\line(1,0){15}}
\put(-15,15){\line(1,-1){15}}
\put(-15,0){\line(1,0){15}}
\put(-15,0){\line(1,1){15}}
\put(-15,0){\line(0,1){15}}
\put(0,0){\line(0,1){15}}

\put(-7.5,23){\line(1,-1){7.5}}
\put(-7.5,23){\line(-1,-1){7.5}}

\put(-7.5,-8){\line(1,1){7.5}}
\put(-7.5,-8){\line(-1,1){7.5}}

\put(-23,7.5){\line(1,1){7.5}}
\put(-23,7.5){\line(1,-1){7.5}}

\put(8,7.5){\line(-1,1){7.5}}
\put(8,7.5){\line(-1,-1){7.5}}

\put(8,-8){\circle*{2}} 
\put(-23,-8){\circle*{2}} 

\put(8,-8){\line(-1,3){7.5}}
\put(8,-8){\line(-3,1){22.5}}

\put(-23,-8){\line(1,3){7.5}}
\put(-23,-8){\line(3,1){22.5}}

\end{picture}
\caption{A 2-vertex connected chordal graph that has no spanning
      2-tree.} 
\label{f2}
\end{figure} 

\begin{rem}
One might be tempted to simplify the proof that
$(1) \Rightarrow (2)$ in Theorem \ref{main1} as follows.
A $k$-tree is either the complete graph on
$k$ vertices, or a graph obtained from a $k$-tree $H$ by
adding a new vertex adjacent to exactly $k$ vertices inducing a $k$-clique in $H$. 
Obviously, a $k$-tree is a chordal graph. 
Now it was shown in \cite{aty10c} that a framework $G(p)$ on $n$ vertices
in general position in
$\Rs^r$, where $G$ is an ($r+1$)-tree admits a positive semidefinite stress matrix
of rank $\rr= n-r-1$. Thus, 
$(1) \Rightarrow (2)$ would follow if every ($r+1$)-vertex connected 
chordal graph has a spanning ($r+1$)-tree. 
Unfortunately, this is not true. Figure \ref{f2} depicts     
a $2$-vertex connected chordal graph which has no spanning $2$-tree.  
\end{rem}

\section{Proof of Theorem \ref{main2}  }
\label{sec4}

Before presenting the proof, we briefly review the connection between Gauss
elimination for sparse symmetric linear systems and chordal graphs. This
connection was discovered by Rose \cite{ros70}. 
Let $A^{(t)}$ denote the matrix $A$ after $t$ steps of Gauss elimination
(one step consists of zeroing out all entries below the pivot).
If no row exchange is needed during the application of
Gauss elimination to an $n \times n$ symmetric 
matrix $A$ of rank $k$, then $A^{(k)}$ is of the form 
\beq \label{Ak} 
A^{(k)} = \left[ \begin{array}{ccccccc}  1 & * & * & * & * & \cdots & * \\  
                                        & 1 & * & * & * & \cdots & * \\
                                        &  & \ddots &* & * & \cdots & * \\
                                        &  &  & 1 & *&  \cdots & *  \\
                                        &  &  &  &  &  & \\
                                        &  &  &  &  &  & \\
                                        &  &  &  &  &  & \\
               \end{array} \right],
\eeq 
where zero entries are left blank and each 
possibly non-zero entry is indicated with an (*).
Note that the entries of the last $n-k$ rows of $A^{(k)}$ are all 0's. 
The assumption that $A$ has generic rank profile, i.e., all the first
$k$ leading principal minors of $A$ are nonzero, ensures that   
Gauss elimination can be applied to $A$ without row exchanges,
and hence $A^{(k)}$ is of the form given  in (\ref{Ak}) \cite{hj85,str09}.

\begin{lem} \label{geA} 
Let $G$ be a chordal graph on $n$ vertices and let
$A$ be an $n \times n$ symmetric matrix of rank $k$ with generic rank profile
such that $a_{ij}=0$ for all $ij: i \neq j$ and $(i,j) \not \in E(G)$. 
Further, assume that $1,2,\ldots,n$ is a PEO of $G$.  
Let $A^{(t)}=(a^{(t)}_{ij})$ denote the matrix $A$ after $t$ steps of Gauss elimination. 
Then,  $a^{(t)}_{ij}=0$ for all $ij:i \neq j$ and $(i,j) \not \in E(G)$, 
and for all $t=1,\ldots,k$. 
\end{lem}

\bpr
Since $A$ has generic rank profile, no row exchanges are needed during
Gauss elimination, i.e., $a_{11} \neq 0$, 
$a^{(1)}_{22} \neq 0, \ldots, a^{(k-1)}_{kk} \neq 0$. 
Therefore, 
\beq \label{gef} 
 a^{(t)}_{ij} = \left\{ \begin{array}{ll} 
            a^{(t-1)}_{ij} & \mbox{if } i < t , \\ 
            a^{(t-1)}_{ij}/a^{(t-1)}_{tt} & \mbox{if } i = t , \\ 
             0 &  \mbox{if }i>t ,j =1,\ldots, t, \\ 
 a^{(t-1)}_{ij}-\frac{a^{(t-1)}_{it} a^{(t-1)}_{tj}}{a^{(t-1)}_{tt}}&
                                 \mbox{if }i>t ,j > t, \\ 
             \end{array} \right.
\eeq
where $A^{(0)}=A$.

Since $G$ is chordal and since $1,2,\ldots,n$ is a PEO of $G$, it follows that
for all $ij: i \neq j$ and $(i,j) \not \in E(G)$, and for all
$s < \min\{i,j\}$, either $(s , i) \not \in E(G)$ or 
$(s , j) \not \in E(G)$ since if  $(s,i) \in E(G)$ and $(s,j) \in E(G)$, then 
$(i,j) \in E(G)$, an obvious contradiction.

The proof is by induction on $t$. 
Let $i\neq j$ and $(i,j) \not \in E(G)$ then $a_{ij}=0$. Thus for  $t=1$ 
we have 
\[
 a^{(1)}_{ij} = \left\{ \begin{array}{ll} 
            a_{ij}/a_{11} & \mbox{if } i = 1 , \\ 
             0 &  \mbox{if }i>1 ,j = 1, \\ 
 a_{ij}-\frac{a_{i1} a_{1j}}{a_{11}}&
                                 \mbox{if }i>1 ,j > 1. \\ 
             \end{array} \right.
\]
Therefore, $a^{(1)}_{ij}=0$  
since $a_{ij}=0$, and since if $i > 1$ and $j > 1$, then either $a_{1i}=0$ or $a_{1j}=0$. 
Therefore, the statement of the lemma is true for $t=1$.  Now assume that the
statement of the lemma is true for $t=m-1$ for some $m:2 \leq m  \leq k$. 
That is, $a^{m-1}_{ij}=0$ for all $ij:i \neq j$ and $(i,j) \not \in E(G)$. 
Then it follows from (\ref{gef}) that 
$a^{(m)}_{ij} = 0$  since  
$a^{(m-1)}_{ij} = 0$, and since if $i > m$ and $j > m$, then 
either $a^{(m-1)}_{im} = 0$ or $a^{(m-1)}_{mj} = 0$.  
Thus the result follows.      
\epr  

\noindent {\bf Proof of Theorem \ref{main2}}

Let $S= \left[ \begin{array}{c} S_1 \\ S_2 \end{array} \right]$
be a stress matrix of $G(p)$ of rank $\rr$ with generic rank profile, 
where $S_1$ is $\rr \times n$ and $S_2$ is $(r+1) \times n$. 
Then
the columns of $S$ span the null space of $\PP$ defined in (\ref{defPe}).
But since $S$ has generic rank profile, it follows that
the first $\rr$ rows of $S$ are linearly independent. Thus
${S_1}^T$ is a Gale matrix of $G(p)$.  

Now 
let $S^{(\rr)}= \left[ \begin{array}{c} {S_1}^{(\rr)} \\ \bz \end{array} \right]$ 
be the matrix $S$ after $\rr$ steps of Gauss elimination, where ${S_1}^{(\rr)}$  
is $\rr \times n$ of rank $\rr$, and $\bz$ is the $(r+1) \times n$ zero matrix; 
and let $Z^T = {S_1}^{(\rr)}$. 
Then $Z$ is a Gale matrix of $G(p)$ since
the rows of $Z^T$ are linear combinations of the rows of ${S_1}^T$.
Furthermore, since $s_{ij}=0$ whenever $i \neq j$ and $(i,j) \not \in E(G)$, 
it follows from Lemma \ref{geA} that
if $i < j$, $i \leq \rr$ and $(i,j) \not \in E(G)$, then $z^T_{ij}=0$. 
Furthermore, $z^T_{ii}=1$ for $i=1,\ldots,\rr$ and $z^T_{ij}=0$ if $i > j$. 
Hence, $Z$ satisfies Property (A).  
Therefore, the result follows from Lemma \ref{Astress}.
\epr

\begin{figure}[t] 
\thicklines 
\setlength{\unitlength}{0.9mm}
\begin{picture}(40,70)(-20,-30) 
\put(-20,0){\circle*{2}} 
\put(-10,-10){\circle*{2}} 
\put(-10,10){\circle*{2}} 
\put(10,-10){\circle*{2}} 
\put(10,10){\circle*{2}} 
\put(20,0){\circle*{2}} 

\put(-20,0){\line(1,1){10}}
\put(-20,0){\line(1,-1){10}}
\put(-20,0){\line(3,-1){30}}
\put(-10,-10){\line(0,1){20}}
\put(-10,-10){\line(1,1){20}}
\put(-10,-10){\line(1,0){20}}
\put(-10,10){\line(1,0){20}}
\put(-10,10){\line(3,-1){30}}
\put(-10,10){\line(1,-1){20}}
\put(10,-10){\line(0,1){20}}
\put(10,-10){\line(1,1){10}}
\put(10,10){\line(1,-1){10}}

\put(-14,-13){$2$}
\put(-24,0){$1$}
\put(-14,11){$3$}
\put(12,-13){$4$}
\put(22,0){$6$}
\put(12,11){$5$}

\end{picture}
\caption{The framework $G(p)$ of Example \ref{exam2}. } 
\label{f4}
\end{figure}
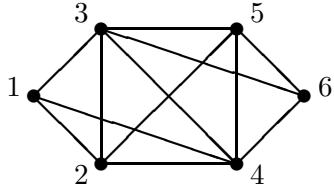 

\begin{exa} \label{exam2}

Consider the chordal framework in general position in $\Rs^2$ shown in 
Figure \ref{f4}, where
$p^1=\left[ \begin{array}{r} -2 \\ 0 \end{array} \right]$,
$p^2=\left[ \begin{array}{r} -1 \\ -1 \end{array} \right]$,
$p^3=\left[ \begin{array}{r} -1 \\ 1 \end{array} \right]$,
$p^4=\left[ \begin{array}{r} 1 \\ -1 \end{array} \right]$,
$p^5=\left[ \begin{array}{r} 1 \\ 1 \end{array} \right]$,
$p^6=\left[ \begin{array}{r} 2 \\ 0 \end{array} \right]$.
Then 
\[
S=\left[ \begin{array}{rrrrrr} 10 & -10 & -5 & 5 & 0 & 0 \\ 
                             -10 & 8 & 7 & -3 & -2 & 0 \\ 
                            -5 & 7 & 1 & -5 & 1 & 1 \\ 
                            5 & -3 & -5 & 1 & 3 & -1 \\ 
                             0 & -2 & 1 & 3 & 0 & -2 \\ 
                            0 & 0 & 1 & -1 & -2 & 2  
     \end{array} \right] 
\]
is an indefinite stress matrix of $G(p)$ of rank 3 with generic rank profile since
the first 3 leading principal minors are nonzero.
Applying Gauss elimination on $S$ leads to
\[
S^{(3)} =\left[ \begin{array}{rrrrrr} 1 & -1 & -0.5 & 0.5 & 0 & 0 \\ 
                              & 1 & -1 & -1 & 1 & 0 \\ 
                             &  & 1 & -1 & -2 & 2 \\ 
                             &  & &  &  &  \\ 
                              &  &  &  &  &  \\ 
                             &  &  &  &  &   
     \end{array} \right]. 
\]
Hence, 
\[
Z =\left[ \begin{array}{rrr} 1 &  & \\ 
                             -1 & 1 &   \\ 
                             -0.5 &-1  & 1 \\ 
                            0.5 &  -1 & -1 \\ 
                             0 & 1 & -2  \\ 
                             0 & 0 &  2 
     \end{array} \right], 
\]
is a Gale matrix of $G(p)$ that satisfies Property (A). 
Furthermore,
\[
ZZ^T =\left[ \begin{array}{rrrrrr} 1 & -1 & -0.5 & 0.5 & 0 & 0 \\ 
                              -1 & 2 & -0.5 & -1.5 & 1 & 0 \\ 
                            -0.5 & -0.5  & 2.25 & -0.25 & -3 & 2 \\ 
                             0.5 &-1.5  & -0.25 & 2.25  & 1  & -2  \\ 
                             0  & 1 & -3  & 1  & 5  & -4  \\ 
                             0 & 0  & 2  & -2  & -4  & 4   
     \end{array} \right] 
\]
is a positive semidefinite stress matrix of $G(p)$ of rank 3.

\end{exa}

\section{Summary and Open Problem}

This paper is a continuation of the study, initiated in \cite{ay10,aty10c},
of the universal rigidity of bar frameworks under the general position assumption. 
In particular, we studied stress matrices of chordal frameworks in general position.
We characterized such frameworks in terms of stress matrices, and   
we showed that for such frameworks universal and global rigidities are equivalent. 
Furthermore, we showed that if a chordal framework on $n$ vertices in general 
position in $r$-dimensions admits
a stress matrix of rank $n-r-1$ with generic rank profile, then it admits a
positive semidefinite stress matrix of rank $n-r-1$. These results suggest
the following problems:

\begin{enumerate}
\item Characterize frameworks $G(p)$ in general position
for which universal and global rigidities are equivalent.  
\item Characterize frameworks $G(p)$ in general position having the property that 
if $G(p)$ has a stress matrix of rank $n-r-1$ with generic rank profile, 
then it has a positive semidefinite stress matrix of rank $n-r-1$.  
\end{enumerate} 

\vspace{0.1in}

\noindent {\bf Acknowledgements:} 
The author would like to thank Prof Leizhen Cai 
for providing the graph in Figure \ref{f2}.

\end{document}